\documentclass[12pt]{article}
\usepackage{amsmath}
\usepackage{cases}
\usepackage{mathrsfs}
\usepackage{multirow}
\usepackage{latexsym,amsmath,amssymb,amsfonts,epsfig,graphicx,cite,psfrag}
\usepackage{eepic,color,colordvi,amscd,amsthm}

\newtheorem{thm}{Theorem}[section]

\newtheorem{lem}[thm]{Lemma}

\newtheorem{prop}[thm]{Proposition}

\numberwithin{equation}{section}

\def\pf{\noindent{\it Proof.} }

\def\qed{\nopagebreak\hfill{\rule{4pt}{7pt}}\medbreak}

\def \C{{\mathcal{C}}}
\def \sg {{\sigma}}
\def \maj {{\mathrm{maj}}}
\def \asc {{\mathrm{asc}}}

\begin{document}
\title{\large\bf Pattern Count on Multiply Restricted Permutations}
\author{Alina F. Y. Zhao\\
\footnotesize School of Mathematical Sciences and Institute of Mathematics \\
\footnotesize Nanjing Normal University, Nanjing 210023, PR China\\
 {\tt alinazhao@njnu.edu.cn}
 }
\date{}
\maketitle

\begin{abstract}
Previous work has studied the pattern count on singly restricted permutations. In this work, we focus on patterns of length $3$ in multiply restricted permutations, especially for double and triple pattern-avoiding permutations. We derive explicit formulae or generating functions for various occurrences of length $3$ patterns on multiply restricted permutations, as well as some combinatorial interpretations for non-trivial pattern relationships.
\end{abstract}

\noindent \textbf{Mathematics Subject Classification:} 05A05, 05A15,
05A19

\section{Introduction}

Let $\sg=\sg_1\sg_2\cdots \sg_n$ be a permutation in the symmetric group $S_n$
written in one-line notation, and  $\sg_i$ is said to be a left to right maximum
(resp. right to left maximum) if $\sg_i>\sg_j$ for all $j<i$ (resp. $j>i$). For a
permutation $q\in S_k$, we say that $\sg$ contains $q$ as a pattern if there exist
$1\leq i_1\leq i_2\leq \cdots \leq i_k \leq n$ such that the entries
$\sg_{i_1}\sg_{i_2}\cdots \sg_{i_k}$ have the same relative order as the entries of
$q$, i.e., $q_j < q_l$ if and only if $\sg_{i_j} < \sg_{i_l}$ whenever $1\leq
j,l\leq k$. We say that $\sg$ avoids $q$ if $\sg$ does not contain $q$ as a pattern.

For a pattern $q$, denote by $S_n(q)$ the set of all permutations in $S_n$ that
avoiding the pattern $q$, and for $R\subseteq S_k$, we denote by $S_n(R)=
\bigcap_{q\in R}S_n(q)$, i.e., the set of permutations in $S_n$ which avoid every
pattern contained in $R$. For two permutations $\sg$ and $q$, we set $f_q(\sg)$ as
the number of occurrences of $q$ in $\sg$ as a pattern, and we further denote the
number of occurrences of $q$ in a permutation set $\Omega$ by
$f_q(\Omega)=\sum_{\sg\in \Omega}f_q(\sg)$.

Recently Cooper \cite{Coo} raised the problem of determining the total number $f_q(S_n(r))$ of
all $q$-patterns in the $r$-avoiding permutations of length $n$. B\'{o}na
\cite{Bona2010} discovered the generating functions of the sequence $f_q(S_n(132))$
for the monotone $q$, and B\'{o}na \cite{Bona} further studied the generating
functions for other length $3$ patterns in the set $S_n(132)$, and showed both
algebraically and bijectively that
\[
f_{231}(S_n(132)) = f_{312}(S_n(132))=f_{213}(S_n(132)).
\]
Rudolph \cite{Rud} proved equipopularity relations between general length $k$ patterns on $132$-avoiding permutations based on the structure of their corresponding binary plane trees. Moreover, Homberger \cite{Hom} also presented exact formulae for each length $3$ pattern in the set $S_n(123)$. Therefore, the singly restricted permutations have been well studied by previous work, whereas it remains open for multiply restricted permutations, e.g., $S_n(123,132)$.

\begin{table*}[t!]
\begin{center}
\begin{tabular}{|p{13mm}|c|c|p{50mm}|}
\hline

\multirow{3}{*}{} & {\tiny $f_{213}(n)=(n-3)2^{n-2}+1 $}

&\multirow{2}{*}{}
&{\tiny $\sum f_{231}(n)x^n=\sum f_{312}(n)x^n=\frac{x^3(1+2x)}{(1-x-x^2)^3}$}\\

\cline{2-2}\cline{4-4}{\tiny $S_n(123,132)$} &
{\tiny$f_{231}(n)=f_{312}(n)=(n^2-5n+8)2^{n-3}-1$} & {\tiny $S_n(123,132,213)$}
&{\tiny $\sum f_{321}(n)x^n= \frac{x^3(1+6x+12x^2+8x^3)}{(1-x-x^2)^4}$}\\
\cline{2-4} & {\tiny $f_{321}(n)=(n^3/3-2n^2+14n/3-5)2^{n-2}+1$} &
&{\tiny $f_{213}(n)=f_{312}(n)={n\choose 3}$} \\
\cline{1-2}\cline{4-4}

\multirow{2}{*}{} & {\tiny $f_{123}(n)=(n-4)2^{n-1}+n+2$} &\multirow{1}{*}{}{\tiny
$S_n(123,132,231)$ } & {\tiny $f_{321}(n)=(n-2){n\choose 3}$}
\\
\cline{2-4}\cline{4-4}{\tiny $S_n(132,213)$} &{\tiny
$f_{231}(n)=f_{312}(n)=(\frac{n^2}{4}-\frac{7n}{4}+4)2^{n}-n-4$} &  &{\tiny
$f_{123}(n)=f_{312}(n)={n+1\choose 4}$}
 \\
\cline{2-2}\cline{4-4}& {\tiny $f_{321}(n)=(\frac{1}{12}n^3-\frac{3}{4}n^2+\frac{38}{12}n-6)2^{n}+n+6$}
&\multirow{1}{*}{}{\tiny $S_n(132,213,231)$} &{\tiny $f_{321}(n)=n(n-2)(n-1)^2/12$}
\\
\hline {\tiny $S_n(132,231)$}& {\tiny
$f_{123}(n)=f_{213}(n)=f_{312}(n)=f_{321}(n)=\frac{2^n}{8}{n\choose 3}$}
& &{\tiny $f_{213}(n)=f_{231}(n)={n\choose 3}$}\\
\cline{1-2}\cline{4-4}

{\tiny $S_n(132,312)$}& {\tiny
$f_{123}(n)=f_{213}(n)=f_{231}(n)=f_{321}(n)=\frac{2^n}{8}{n\choose 3}$}
& \multirow{1}{*}{}{\tiny $S_n(123,132,312)$} &{\tiny $f_{321}(n)=(n-2){n\choose 3}$}\\
\hline

\multirow{2}{*}{}
& {\tiny $f_{213}(n)=f_{231}(n)=f_{312}(n)={n+2\choose 5}$}&&{\tiny $f_{132}(n)=f_{213}(n)={n+1\choose 4}$}\\
\cline{2-2}\cline{4-4}{\tiny $S_n(132,321)$}&  {\tiny
$f_{123}(n)=\frac{7n^5}{120}-\frac{n^4}{3}+\frac{17n^3}{24}-\frac{2n^2}{3}+\frac{7}{30}$}
&\multirow{1}{*}{}{\tiny $S_n(123,231,312)$}&{\tiny $f_{321}(n)=\frac{1}{12}n(n-2)(n-1)^2$} \\
\hline
\end{tabular}
\end{center}
\caption{The pattern count on doubly and triply restricted permutations.}\label{tab:sum}
\end{table*}

In this paper, we are interested in pattern count on multiply restricted permutations $S_n(R)$ for $R\subset S_3$, especially for double and triple restrictions. We derive explicit formulae or generating functions for the occurrences of each length $3$ pattern in multiply restricted permutations, and the detailed results are summarized in Table~\ref{tab:sum}. Also, we present some combinatorial interpretations for non-trivial pattern relationships. It is trivial to consider the restricted permutations of higher multiplicity since there are only finite permutations, as shown in \cite{Sim}. Therefore, this work presents a complete study on length $3$ patterns of multiply restricted permutations.

\section{Doubly Restricted Permutations}
For $\sg \in S_n$, the following three operations are very useful in pattern avoiding enumeration. The complement of $\sg$ is given by $\sg^c=(n+1-\sg_1)(n+1-\sg_2)\cdots (n+1-\sg_n)$, its reverse is defined as $\sg^r=\sg_n\cdots \sg_2\sg_1$ and the inverse $\sg^{-1}$ is the group theoretic inverse permutation. For any set of permutations $R$, let $R^c$ be the set obtained by complementing each element of $R$, and the set $R^r$ and $R^{-1}$ are defined analogously.
\begin{lem}\label{oper}
Let $R\subseteq S_k$ be any set of permutations in $S_k$, and $\sg \in S_n$, we have
\[ \sg\in S_n(R)\Leftrightarrow \sg^c\in S_n(R^c)\Leftrightarrow \sg^r\in S_n(R^r)\Leftrightarrow \sg^{-1}\in S_n(R^{-1}).\]
\end{lem}

From Lemma~\ref{oper} and the known results on $S_n(123)$ and $S_n(132)$, we can
obtain each length $3$ pattern count in the singly restricted permutations
$S_n(r)$ for $r=213,231,312$ and $321$. In this section, we focus on the number of
length $3$ patterns in all doubly restricted permutations.

A composition of $n$ is an expression of $n$ as an ordered sum of positive integers,
and we say that $c$ has $k$ parts or $c$ is a $k$-composition if there are exactly
$k$ summands appeared in a composition $c$. Denote by $\mathcal{C}_n$ and
$\mathcal{C}_{n,k}$ the set of all compositions of $n$ and the set of $k$-compositions
of $n$, respectively. It is known that $|\mathcal{C}_n|=2^{n-1}$ and
$|\mathcal{C}_{n,k}|= {n-1\choose k-1}$ for $n\geq 1$ with $|\mathcal{C}_0|=1$. For
more details, see \cite{Sta}.

We begin with some enumerative results on compositions.
\begin{lem}\label{com1}
For $n\geq 1$, we have
\begin{eqnarray}
a(n)&:=&\sum_{k\geq 1 \atop c_1+\cdots+c_{k-1}+c_k=n}c_k=2^n-1,\\
b(n)&:=&\sum_{k\geq 1 \atop c_1+\cdots+c_{k-1}+c_k=n}c_k(c_k-1)=2^{n+1}-2n-2,
\end{eqnarray}
where the sum takes over all compositions of $n$.
\end{lem}
\pf For $c_k=m$ , we can regard $c_1+c_2+\cdots+c_{k-1}$ as a composition of $n-m$. It is easy to get that the number of compositions of $n-m$ is $2^{n-m-1}$ for $1\leq m\leq n-1$ and the number of empty compositions is one. This follows that
\begin{align*}
a(n)&=n+\sum_{m=1}^{n-1}m2^{n-m-1},\\
b(n)&=n(n-1)+\sum_{m=1}^{n-1}m(m-1)2^{n-m-1}.
\end{align*}
Let $g(x)=\sum_{i=0}^{n-1}x^{i}=\frac{1-x^n}{1-x}$, and we have
\begin{eqnarray*}
&g'(x)=\sum_{i=1}^{n-1}ix^{i-1}=\frac{(n-1)x^n-nx^{n-1}+1}{(1-x)^2}, \\
&g''(x)=\sum_{i=1}^{n-1}i(i-1)x^{i-2}=\frac{(3n-n^2-2)x^n+(2n^2-4n)x^{n-1}+(n-n^2)x^{n-2}+2}{(1-x)^3}.&
\end{eqnarray*}
By setting $x=1/2$ in $g'(x)$ and $g''(x)$, we get
\begin{eqnarray*}
g'(1/2)&=&2^2((n-1)2^{-n}-n2^{-n+1}+1),\\
g''(1/2)&=&2^3\left[(3n-n^2-2)2^{-n}+(2n^2-4n)2^{-n+1}+(n-n^2)2^{-n+2}+2\right].
\end{eqnarray*}
Observing $a(n)=2^{n-2}g'(1/2)+n$ and $b(n)=2^{n-3}g''(1/2)+n(n-1)$, this lemma follows as desired.\qed

\begin{lem}\label{com2}
For $n\geq 1$, we have
\begin{eqnarray}
c(n)&:=&\sum_{k\geq 1 \atop c_1+\cdots+c_{k-1}+c_k=n}k=(n+1)2^{n-2},\\
d(n)&:=&\sum_{k\geq 1 \atop c_1+\cdots+c_{k-1}+c_k=n}k(k-1)=(n^2+n-2)2^{n-3},
\end{eqnarray}
where the sum takes over all compositions of $n$.
\end{lem}
\pf Since the number of compositions of $n$ with $k$ parts is
${n-1\choose k-1}$, we have
\[c(n)=\sum_{k=1}^{n}k{n-1\choose k-1} \text{ and } d(n)=\sum_{k=1}^{n}k(k-1){n-1\choose k-1}.\]
Let $h(x)=x\sum_{i=1}^{n}{n-1\choose i-1}x^{i-1}=x(1+x)^{n-1}$. Then
\begin{eqnarray*}
h'(x) &=&\sum_{i=1}^n i{n-1\choose i-1} x^{i-1}=(nx+1)(1+x)^{n-2},  \\
h''(x)&=&\sum_{i=1}^n i(i-1){n-1\choose i-1}x^{i-2} =\big(n^2x+n(2-x)-2\big)(1+x)^{n-3}.
\end{eqnarray*}
We complete the proof by putting $x=1$ in the above formulae.\qed

Based on Lemma~\ref{oper}, Simion and Schmidt \cite{Sim} showed that the pairs of patterns among the total ${6\choose 2}=15$ cases fall into the following $6$ classes.
\begin{prop}[\cite{Sim}]
For every symmetric group $S_n$,
{\small
\begin{enumerate}
\item $|S_n(123,132)|= |S_n(123,213)|=|S_n(231,321)|= |S_n(312,321)|= 2^{n-1}$;
\item $|S_n(132,213)|=|S_n(231,312)|=2^{n-1}$;
\item $|S_n(132,231)|=|S_n(213,312)|=2^{n-1}$;
\item $|S_n(132,312)|=|S_n(213,231)|=2^{n-1}$;
\item $|S_n(132,321)|=|S_n(123,231)|=|S_n(123,312)|=|S_n(213,321)|={n\choose 2}+1$;
\item $|S_n(123,321)|=0$ for $n\geq 5$.
\end{enumerate}}
\end{prop}

Therefore, it is sufficient to consider the pattern count of the first set for each class in the subsequent sections, and the other sets can be obtained by taking the complement or reverse or inverse of the known results.

\subsection{Pattern Count on $(123,132)$-Avoiding Permutations}
We first present a bijection between $S_n(123,132)$ and $\mathcal{C}_n$ as follows:
\begin{lem}\label{ta}
There is a bijection $\varphi_1$ between the sets $S_n(123,132)$ and $\mathcal{C}_n$.
\end{lem}
\pf For any given $\sg\in S_n(123,132)$, let $\sg_{i_1},\sg_{i_2},\ldots,\sg_{i_k}$ be the $k$ right to left maxima with $i_1<i_2<\cdots<i_k$, which yields that $c=i_1+(i_2-i_1)+\cdots+(i_{k-1}-i_{k-2})+(i_{k}-i_{k-1})$ is a composition of $n$ since $i_k=n$. On the converse, let $m_i=n-(c_1+\cdots+c_{i-1})$ for any given $n=c_1+c_2+\cdots+c_k\in \C_n$, and we set $\tau_i=m_i-1,m_i-2,\ldots,m_i-c_i+1,m_i$ for $1\leq i \leq k$. Therefore, $\sg=\tau_1,\tau_2,\ldots,\tau_k \in S_n(123,132)$ is as desired. For example, we have $\sg=8\,9\,7\,5\,4\,3\,6\,1\,2$ for the composition $9=2+1+4+2$. \qed

For a pattern $q$, we denote by $f_{q}(n):=\sum_{\sigma\in
S_n(123,132)}f_q(\sigma)$, i.e., the number of occurrences of the pattern $q$ in $S_n(123,132)$. For simplicity, we will use this notation in subsequent sections when the set in question is unambiguous. For convenience, we denote by
$\tau_i>\tau_j$ if all the elements in the subsequence $\tau_i$ are larger than all
the elements in subsequence $\tau_j$. Based on Lemma~\ref{ta}, we have
\begin{prop}\label{pa}
For $n\geq3$,
\begin{eqnarray}
f_{213}(n)&=&\sum_{k\geq1 \atop c_1+c_2+\cdots+c_k=n}\sum_{i=1}^k{c_i-1\choose 2},\label{pa1}\\
f_{231}(n)&=&\sum_{k\geq1 \atop c_1+c_2+\cdots+c_k=n}\sum_{i=1}^{k-1}\sum_{j=
i+1}^kc_j(c_i-1) .\label{pa2}
\end{eqnarray}
\end{prop}
\pf For each permutation $\sg\in S_n(123,132)$ with $\varphi_1(\sg)=c_1+ c_2+
\cdots+c_k$, we can rewrite $\sg$ as $\sg=\tau_1,\tau_2,\ldots,\tau_k$ from
Lemma~\ref{ta}. For $j>i$, since $\tau_i>\tau_{j}$, and the elements except the last one are decreasing in $\tau_i$, the pattern $213$ can only occur in every subsequence $\tau_i$. Thus, we have ${c_i-1 \choose 2}$
choices to choose two elements in $\tau_i$ to play the role of $21$, and the last element of $\tau_i$ plays the role of $3$, summing up all
the number of $213$-patterns in subsequences $\tau_1,\tau_2,\ldots,\tau_k$ gives the
formula \eqref{pa1}.

For the pattern $231$, we have $c_i-1$ choices in the subsequence $\tau_i$ to choose one element playing the role of $2$ and one choice (always the last element of $\tau_i$) for $3$, and then, we have $c_{i+1}+\cdots+c_k$ choices to choose one element for the role of $1$ since all the elements after $\tau_i$ are smaller than those in $\tau_i$. Summing up all the number of $231$-patterns according to the position of $3$ gives the formula \eqref{pa2}. \qed

Based on the previous analysis, we now present our first main results of the
explicit formulae for pattern count in the set $S_n(123,132)$.

\begin{thm}\label{ta1}
For $n\geq 3$, in the set $S_n(123,132)$, we have
\begin{eqnarray}
f_{213}(n)&=&(n-3)2^{n-2}+1, \label{fta1} \\
f_{231}(n)&=&f_{312}(n)=(n^2-5n+8)2^{n-3}-1,\label{fta2}   \\
f_{321}(n)&=&(n^3/3-2n^2+14n/3-5)2^{n-2}+1.\label{fta4}
\end{eqnarray}
\end{thm}

\pf From $S_3(123,132)= \{213,231,312,321\}$, it is obvious that
\[
f_{213}(3)=f_{231}(3)=1.
\]
To prove formula~\eqref{fta1}, from Prop~\ref{pa}, we observe that, for $n\geq3$
\[
f_{213}(n+1)=\sum_{k\geq 1, c_k=1 \atop c_1+c_2+\cdots+c_k=n+1}\sum_{i=1}^k {c_i-1\choose2}+\sum_{k\geq 1, c_k\geq2 \atop c_1+c_2+\cdots+c_k=n+1}\sum_{i=1}^k {c_i-1\choose2}.
\]
If $c_k=1$, then $k\geq2$, and we further have
\[
\sum_{k\geq 1, c_k=1 \atop c_1+c_2+\cdots+c_k=n+1}\sum_{i=1}^k
{c_i-1\choose2}=\sum_{k-1\geq 1\atop c_1+c_2+\cdots+c_{k-1}=n}\sum_{i=1}^{k-1}
{c_i-1\choose 2}=f_{213}(n).
\]
If $c_k\geq 2$, then we set $c_k=1+r_k$, and from Lemma~\ref{com1}, it holds that
\begin{align*}
\sum_{k\geq 1, c_k\geq2 \atop c_1+c_2+\cdots+c_k=n+1}\sum_{i=1}^k {c_i-1\choose2}&=\sum_{k\geq 1 \atop c_1+\cdots+c_{k-1}+r_k=n}\left[\sum_{i=1}^{k-1}{c_i-1\choose 2}+{r_k-1\choose 2}+(r_k-1)\right]\\
&=f_{213}(n)+\sum_{k\geq 1\atop c_1+\cdots+c_{k-1}+r_k=n}(r_k-1) \\ &=f_{213}(n)+a(n)-\sum_{k\geq 1\atop c_1+\cdots+c_{k-1}+r_k=n}1
=f_{213}(n)+a(n)-2^{n-1}.
\end{align*}
Combining the above two cases, we have
\[f_{213}(n+1)=2f_{213}(n)+2^{n-1}-1.\]
This proves formula~\eqref{fta1} by solving the recurrence with initial value $f_{213}(3)=1$.

To prove formula~\eqref{fta2}, first observe that from Lemma~\ref{oper}, $\sg \in
S_n(123,132)\Leftrightarrow \sg^{-1}\in S_n(123,132)$, this implies the first
equality of formula~\eqref{fta2} directly since $231^{-1}=312$. While for the second
equality of formula~\eqref{fta2}, by Prop~\ref{pa}, we have
\[
f_{231}(n+1)=\sum_{k\geq1,c_k=1 \atop
c_1+c_2+\cdots+c_k=n+1}\sum_{i=1}^{k-1}\sum_{j= i+1}^k
c_j(c_i-1)+\sum_{k\geq1,c_k\geq2 \atop
c_1+c_2+\cdots+c_k=n+1}\sum_{i=1}^{k-1}\sum_{j= i+1}^k c_j(c_i-1).
\]
If $c_k=1$, then $k\geq 2$, and from Lemma~\ref{com2} we have
\begin{align*}
\sum_{k\geq1,c_k=1 \atop c_1+c_2+\cdots+c_k=n+1}\sum_{i=1}^{k-1}\sum_{j= i+1}^k c_j(c_i-1)&=\sum_{k-1\geq 1\atop c_1+\cdots+c_{k-1}=n}\sum_{i=1}^{k-2}\sum_{j= i+1}^{k-1} c_j(c_i-1) +\sum_{k-1\geq 1\atop c_1+\cdots+c_{k-1}=n}\sum_{i=1}^{k-1}(c_i-1)\\
&= f_{231}(n)+\sum_{k-1\geq 1 \atop c_1+\cdots+c_{k-1}=n}n-(k-1)\\
&= f_{231}(n)-c(n)+n2^{n-1}.
\end{align*}
If $c_k\geq 2$, then we set $c'_k=c_k-1$ and $c'_i=c_i$ for $1\leq i\leq k-1$. This holds that
\begin{align*}
\sum_{k\geq1,c_k\geq2 \atop c_1+\cdots+c_k=n+1}\sum_{i=1}^{k-1}\sum_{j= i+1}^k c_j(c_i-1)&=\sum_{k\geq1 \atop c'_1+\cdots+c'_k=n}\sum_{i=1}^{k-1}\sum_{j= i+1}^k c'_j(c'_i-1)+\sum_{k\geq1 \atop c'_1+\cdots+c'_k=n}\sum_{i=1}^{k-1}(c'_i-1)\\
&= f_{231}(n)+\sum_{k\geq 1\atop c'_1+\cdots+c'_k=n}(n-c'_k-k+1)\\
&=f_{231}(n)-a(n)-c(n)+(n+1)2^{n-1},
\end{align*}
where the last equality holds from Lemma~\ref{com1} and Lemma~\ref{com2}. Therefore, after simplification, we have
\[
f_{231}(n+1)=2f_{231}(n)+(n-2)2^{n-1}+1,
\]
which proves the second equality of formula~\eqref{fta2} by using $f_{231}(3)=1$.

Note that the total number of all length $3$ patterns in a permutation $\sg \in S_n$ is ${n\choose 3}$, for the set $S_n(123,132)$, this gives the relation
\[
f_{213}(n)+2f_{231}(n)+f_{321}(n)={n\choose 3}2^{n-1}.
\]
Thus formula~\eqref{fta4} is a direct computation of the above equation. This completes the proof.\qed

The first few values of $f_q(S_n(123,132))$ for $q$ of length $3$ are shown below.
\begin{center}
\begin{tabular}{|c|c|c|c|c|c|c||c|c|c|c|c|c|c|}
\hline   $n$ &$f_{123}$&$f_{132}$&$f_{213}$ &$f_{231}$&$f_{312}$  &$f_{321}$
       & $n$ &$f_{123}$&$f_{132}$&$f_{213}$ &$f_{231}$&$f_{312}$  &$f_{321}$\\
\hline   $3$ &$0$      &$0$      &$1$       &$1$      &$1$&$1$ &$6$ &$0$      &$0$       &$49$     &$111$    &$111$     &$369$\\
\hline   $4$ &$0$      &$0$       &$5$       &$7$      &$7$       &$13$ &$7$ &$0$      &$0$       &$129$     &$351$    &$351$     &$1409$\\
\hline   $5$ &$0$      &$0$       &$17$      &$31$     &$31$      &$81$ &$8$ &$0$      &$0$       &$321$     &$1023$    &$1023$     &$4801$\\
\hline
\end{tabular}\label{2a}
\end{center}

\subsection{Pattern Count on $(132,213)$-Avoiding Permutations}
We begin with the following correspondence between $(132,213)$-avoiding
permutations and compositions of $n$.
\begin{lem}\label{tb}
There is a bijection $\varphi_2$ between the sets $S_n(132,213)$ and
$\mathcal{C}_n$.
\end{lem}
\pf Given $\sg\in S_n(132,213)$, let $\sg_{i_1},\sg_{i_2},\ldots, \sg_{i_k}$ be the
$k$ right to left maxima with $i_1<i_2<\cdots<i_k$. This follows that
$c=i_1+(i_2-i_1)+\cdots+(i_{k-1}-i_{k-2})+(i_{k}-i_{k-1})$ is a composition of $n$
since $i_k=n$. On the converse, given $n=c_1+c_2+\cdots+c_k\in \C_n$, let
$m_i=n-(c_1+\cdots+c_{i-1})$ and $\tau_i=m_i-c_i+1,m_i-c_i+2,\ldots,m_i-1,m_i$ for
$1\leq i \leq k$. We set $\sg=\tau_1,\tau_2,\ldots,\tau_k$, and it is easy to check
that $\sg \in S_n(132,213)$. For example, for the composition $9=3+3+1+2$, we get
$\sg=7\,8\,9\,4\,5\,6\,3\,1\,2$. \qed

Based on the above lemma, we have
\begin{prop}\label{pb}
For $n\geq3$,
\begin{align}
f_{123}(n)&=\sum_{k\geq 1 \atop c_1+c_2+\cdots+c_k=n} \sum_{i=1}^k{c_i\choose 3},\label{pb1}\\
f_{231}(n)&=\sum_{k\geq 1 \atop c_1+c_2+\cdots+c_k=n}\sum_{i=1}^{k-1}
\sum_{j=i+1}^kc_j {c_i\choose 2}\label{pb2}.
\end{align}
\end{prop}
\pf For a permutation $\sg\in S_n(132,213)$ with
$\varphi_2(\sg)=c_1+c_2+\cdots+c_k$, we rewrite $\sg$ as
$\sg=\tau_1,\tau_2,\ldots,\tau_k$ by Lemma~\ref{tb}, and we see that the pattern
$123$ can only occur in every subsequence $\tau_i$ since $\tau_i>\tau_{j}$ for $j>i$
and the elements in $\tau_i$ are increasing. Thus, we have ${c_i \choose 3}$ choices
to choose three elements in $\tau_i$ to play the role of $123$, and formula
\eqref{pb1} follows by summing all $123$-patterns in subsequences
$\tau_1,\tau_2,\ldots,\tau_k$.

For pattern $231$, we have ${c_i \choose 2}$ choices in
the subsequence $\tau_i$ to choose two elements to play the role of
$23$, after this we have $c_{i+1}+\cdots+c_k$ choices to choose one
element in $\tau_{i+1},\ldots,\tau_k$ for the role of $1$ since
$\tau_j<\tau_i$ for all $j>i$. Summing up all the number of
$231$-patterns according to the positions of $23$ gives the formula
\eqref{pb2}. \qed

\begin{thm}\label{tb1}
For $n\geq 3$, in the set $S_n(132,213)$, we have
\begin{align}
f_{123}(n)&=(n-4)2^{n-1}+n+2,\label{ftb1}\\
f_{231}(n)&=f_{312}(n)=(n^2-7n+16)2^{n-2}-n-4,\label{ftb2}\\
f_{321}(n)&=(n^3/3-3n^2+38n/3-24)2^{n-2}+n+6.\label{ftb3}
\end{align}
\end{thm}

\pf From Prop~\ref{pb}, we have
\[
f_{123}(n+1)=\sum_{k\geq 1, c_k=1 \atop c_1+c_2+\cdots+c_k=n+1}\sum_{i=1}^k
{c_i\choose 3}+\sum_{k\geq 1, c_k\geq2 \atop c_1+c_2+\cdots+c_k=n+1}\sum_{i=1}^k
{c_i\choose 3}.
\]
If $c_k=1$, then $k\geq 2$, and we have
\begin{align*}
\sum_{k\geq 1 ,c_k=1 \atop c_1+c_2+\cdots+c_{k}=n+1}\sum_{i=1}^{k}
{c_i\choose 3}=\sum_{k-1\geq 1 \atop c_1+c_2+\cdots+c_{k-1}=n}\sum_{i=1}^{k-1}
{c_i\choose 3}=f_{123}(n).
\end{align*}
If $c_k\geq 2$, then we set $c_k=1+r_k$, where $r_k\geq 1$, and this follows
\begin{align*}
\sum_{k\geq 1, c_k\geq2 \atop c_1+c_2+\cdots+c_k=n+1}\sum_{i=1}^k {c_i\choose
3}&=\sum_{k\geq 1 \atop c_1+c_2+\cdots+c_{k-1}+r_k=n}
\left[\sum_{i=1}^{k-1}{c_i\choose 3}+{r_k\choose 3}+\frac{r_k(r_k-1)}{2}\right]\\
&=f_{123}(n)+b(n)/2.
\end{align*}
Combining the above two cases, we get that
\[f_{123}(n+1)=2f_{123}(n)+2^n-n-1,\]
and formula~\eqref{ftb1} holds by solving the recurrence with initial value $f_{123}(3)=1$.

From Lemma~\ref{oper}, we see that $\sg \in S_n(132,213)\Leftrightarrow \sg^{-1}\in
S_n(132,213)$, and this follows that $f_{231}(n)=f_{312}(n)$ from  $231^{-1}=312$.

To calculate $f_{231}(n)$, we have by using Prop~\ref{pb} again
\[f_{231}(n+1)=\sum_{k\geq 1,c_k=1 \atop c_1+c_2+\cdots+c_k=n+1}\sum_{i=1}^{k-1}
\sum_{j=i+1}^kc_j {c_i\choose 2}+\sum_{k\geq 1,c_k\geq 2 \atop c_1+c_2+\cdots+c_k=n+1}\sum_{i=1}^{k-1}
\sum_{j=i+1}^kc_j {c_i\choose 2}.
\]
If $c_k=1$, then $k\geq 2$, and
\begin{align*}
\sum_{k\geq 1,c_k=1 \atop c_1+c_2+\cdots+c_k=n+1}\sum_{i=1}^{k-1}
\sum_{j=i+1}^kc_j{c_i\choose 2} &=\sum_{k-1\geq 1 \atop
c_1+c_2+\cdots+c_{k-1}=n}\sum_{i=1}^{k-1}
{c_i\choose 2}\left[\sum_{j=i+1}^{k-1}c_j+1\right]\\
&=f_{231}(n)+\alpha(n),
\end{align*}
where $\alpha(n)=\sum \limits_{k\geq 1\atop c_1+c_2+\cdots+c_k=n}\sum_{i=1}^{k}
{c_i\choose 2}$. We further have
\begin{align*}
\alpha(n)&=\sum_{k\geq 1 \atop c_1+\cdots+c_{k}=n}\sum_{i=1}^{k}\left[{c_i-1\choose 2}+c_i-1\right]
=f_{213}(S_n(123,132))+\sum_{k\geq 1 \atop c_1+\cdots+c_{k}=n}(n-k)\\
&=f_{213}(S_n(123,132))-c(n)+n2^{n-1},
\end{align*}
in which we use the derived formula $f_{213}(S_n(123,132))=\sum\limits_{k\geq 1 \atop
c_1+\cdots+c_{k}=n}\sum_{i=1}^{k}{c_i-1\choose 2}$.\\
If $c_k\geq 2$, then we set $c'_k=c_k-1$ and $c'_i=c_i$ for $1\leq i\leq k-1$. This
holds that
\begin{align*}
\sum_{k\geq 1,c_k\geq 2 \atop c_1+\cdots+c_k=n+1}
\sum_{i=1}^{k-1}\sum_{j=i+1}^kc_j{c_i\choose 2} &= \sum_{k\geq 1 \atop
c'_1+\cdots+c'_k=n} \sum_{i=1}^{k-1}\sum_{j=i+1}^kc'_j{c'_i\choose 2}+ \sum_{k\geq 1
\atop c'_1+\cdots+c'_k=n}
\sum_{i=1}^{k-1}{c'_i\choose 2}\\
&=f_{231}(n)+\beta(n),
\end{align*}
where $\beta(n)=\sum\limits_{k\geq 1 \atop c_1+\cdots+c_k=n}
\sum_{i=1}^{k-1}{c_i\choose 2}$. Further, we can rewrite $\beta(n)$ as
\begin{align*}
\beta(n)&=\sum_{k\geq 1 \atop c_1+\cdots+c_k=n}
\sum_{i=1}^{k}{c_i\choose 2}-\sum_{k\geq 1 \atop c_1+\cdots+c_k=n}
\frac{c_k (c_{k}-1)}{2}\\
&=\alpha(n)-b(n)/2=f_{213}(S_n(123,132))-c(n)-b(n)/2+n2^{n-1}.
\end{align*}
Substituting the known formulae for $f_{213}(S_n(123,132))$, $c(n)$ and $b(n)$, we
get that
\[
f_{231}(n+1)=2f_{231}(n)+(2n-6)2^{n-1}+n+3,
\]
and formula~\eqref{ftb3} holds by solving this recurrence with initial condition
$f_{213}(3)=1$.

Finally, formula \eqref{ftb3} follows from $f_{123}(n)+2f_{231}(n)+f_{321}(n)={n\choose 3}2^{n-1}$, and this completes the proof. \qed

The first few values of $f_q(S_n(132,213))$ for $q$ of length $3$ are shown below.
\begin{center}
\begin{tabular}{|c|c|c|c|c|c|c||c|c|c|c|c|c|c|}
\hline $n$ &$f_{123}$&$f_{132}$&$f_{213}$ &$f_{231}$&$f_{312}$  &$f_{321}$ &$n$ &$f_{123}$&$f_{132}$&$f_{213}$ &$f_{231}$&$f_{312}$  &$f_{321}$\\
\hline   $3$ &$1$      &$0$      &$0$       &$1$      &$1$&$1$ &$6$ &$72$      &$0$       &$0$     &$150$    &$150$     &$268$\\
\hline   $4$ &$6$      &$0$       &$0$       &$8$      &$8$       &$10$ &$7$ &$201$      &$0$       &$0$     &$501$    &$501$     &$1037$\\
\hline   $5$ &$23$      &$0$       &$0$      &$39$     &$39$      &$59$ &$8$ &$522$      &$0$       &$0$     &$1524$    &$1524$     &$3598$\\
\hline
\end{tabular}\label{2b}
\end{center}

\subsection{Pattern Count on $(132,231)$-Avoiding Permutations}
\begin{thm}\label{tc1}
For $n\geq 3$, in the set $S_n(132,231)$, we have
\begin{align}\label{ftc1}
f_{123}(n)=f_{213}(n)=f_{312}(n)=f_{321}(n)={n\choose 3}2^{n-3}.
\end{align}
\end{thm}

\pf For each $\sg \in S_n(132,231)$, we observe that $n$ must lie in the beginning
or the end of $\sg$, and $n-1$ must lie in the beginning or the end of $\sg
\backslash \{n\}$, ..., and so on. Here $\sg \backslash \{n\}$ denotes the sequence
obtained from $\sg$ by deleting element $n$. Based on such observation, suppose
$abc$ is a length $3$ pattern in $S_n(132,231)$, and set $[n] \backslash
\{a,b,c\}:=\{r_1>r_2>\cdots>r_{n-4}>r_{n-3}\}$. We can construct a permutation in
the set $S_n(132,231)$ which contains an $abc$ pattern as follows:

Start with the subsequence $\sg^0:=abc$, and for $i$ from $1$ to $n-3$, $\sg^i$ is
obtained from $\sg^{i-1}$ by inserting $r_i$ into it.
\begin{itemize}
\item if there are at least two elements in $\sg^{i-1}$ smaller than $r_i$, then choose the two elements $A$ and $B$ such that $A$ is the leftmost one and $B$ is the rightmost one. We put $r_i$ immediately to the left of $A$ or immediately to the right of $B$;
\item if there is only one element $A$ in $\sg^{i-1}$ such that $A<r_i$, then we can put $r_i$ immediately to the left or to the right of $A$;
\item if all the elements in $\sg^{i-1}$ are larger than $r_i$, then choose $A$ the least one, and put $r_i$ immediately to the left or to the right of $A$.
\end{itemize}

Finally, we set $\sg:=\sg^{n-3}$ and $\sg\in S_n(132,231)$ from the above
construction. Moreover, there are ${n\choose 3}$ choices to choose $abc$, and the
number of permutations having $abc$ as a pattern is $2^{n-3}$ since each $r_i$ has
$2$ choices. This completes the proof.\qed

Here we give an illustration of constructing a permutation in $S_8(132,231)$ which
contains the pattern $abc=256$. Set $\sg^0:=256$, we may have $\sg^1=8\,2\,5\,6$,
$\sg^2=8\,7\,2\,5\,6$, $\sg^3=8\,7\,2\,4\,5\,6$, $\sg^4=8\,7\,3\,2\,4\,5\,6$,
$\sg:=\sg^5=8\,7\,3\,2\,1\,4\,5\,6$.

We could also provide combinatorial proofs for the phenomenon
$f_{123}(n)=f_{213}(n)=f_{312}(n)=f_{321}(n)$. From Lemma~\ref{oper}, we have $\sg
\in S_n(132,231)\Leftrightarrow \sg^{r}\in S_n(132,231)$, and this follows
$f_{123}(n)=f_{321}(n)$ and $f_{213}(n)=f_{312}(n)$ from $123^{r}=321$ and
$213^{r}=312$, respectively. It remains to give a bijection for
$f_{213}(n)=f_{123}(n)$, and our following construction is motivated from B\'{o}na
\cite{Bona}.

A binary plane tree is a rooted unlabeled tree in which each vertex has at most two children, and each child is a left child or a right child of its parent. For each $\sg\in S_n(132)$, we construct a binary plane tree $T(\sg)$ as follows: the root of $T(\sg)$
corresponds to the entry $n$ of $\sg$, the left subtree of the root
corresponds to the string of entries of $\sg$ on the left of $n$,
and the right subtree of the root corresponds to the string of
entries of $\sg$ on the right of $n$. Both subtrees are constructed
recursively by the same rule. For more details, see
\cite{Bona2011,Bona,Rud}.

A left descendant (resp. right descendant) of a vertex $x$ in a
binary plane tree is a vertex in the left (resp. right) subtree of
$x$. The left (resp. right) subtree of $x$ does not contain $x$
itself. Similarly, an ascendant of a vertex $x$ in a binary plane
tree is a vertex whose subtree contains $x$. Given a tree $T$ and a
vertex $v \in T$, let $T_v$ be the subtree of $T$ with $v$ as the
root. Let $R$ be an occurrence of the pattern $123$ in $\sg \in
S_n(132)$, and let $R_1,R_2,R_3$ be the three vertices of $T(\sg)$
that correspond to $R$, going left to right. Then, $R_1$ is a left
descendant of $R_2$, and $R_2$ is a left descendant of $R_3$.

From the above correspondence, we see that for $\sg \in S_n(132,231)$, $T(\sg)$ is a
binary plane tree on $n$ vertices such that each vertex has at most one child from
the forbiddance of the pattern $231$. For simplicity, denote by $\mathcal{T}_n$ the
set of such binary plane trees on $n$ vertices. Let $Q$ be an occurrence of the
pattern $213$ in $\sg \in S_n(132,231)$, and let $Q_2,Q_1,Q_3$ be the three vertices
of $T(\sg)$ that correspond to $Q$, going left to right. From the characterization
of trees in $\mathcal{T}_n$,  $Q_2$ is a left descendant of $Q_3$, and $Q_1$ is a
right descendant of $Q_2$.

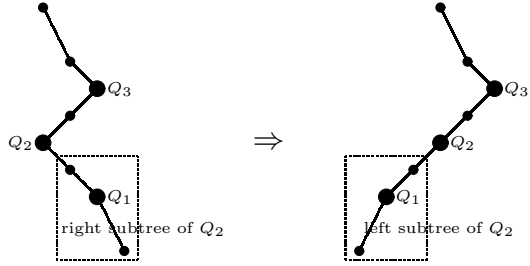
\begin{figure}[h,t]
\setlength{\unitlength}{0.6mm}%
\begin{center}
\begin{picture}(50,55)
\put(0,24){\circle*{3.5}}\put(6,18){\circle*{2}}
\put(12,12){\circle*{3.5}}\put(18,0){\circle*{2}}
\put(6,30){\circle*{2}} \put(12,36){\circle*{3.5}}
\put(6,42){\circle*{2}} \put(0,54){\circle*{2}}
\qbezier[1000](0,24)(3,21)(6,18) \qbezier[1000](6,18)(9,15)(12,12)
\qbezier[1000](12,12)(15,6)(18,0) \qbezier[1000](0,24)(3,27)(6,30)
\qbezier[1000](6,30)(9,33)(12,36) \qbezier[1000](12,36)(9,39)(6,42)
\qbezier[1000](0,54)(3,48)(6,42) 
\put(17,12){\makebox(0,0){\tiny $Q_1$}}
\put(-5,24){\makebox(0,0){\tiny $Q_2$}}
\put(17,36){\makebox(0,0){\tiny $Q_3$}}

\put(50,24){\makebox(0,0){$\Rightarrow$}}
\put(70,0){\circle*{2}}\put(76,12){\circle*{3.5}}
\put(82,18){\circle*{2}} \put(88,24){\circle*{3.5}}
\put(94,30){\circle*{2}} \put(100,36){\circle*{3.5}}
\put(94,42){\circle*{2}} \put(88,54){\circle*{2}}
\qbezier[1000](70,0)(73,6)(76,12)
\qbezier[1000](76,12)(79,15)(82,18)\qbezier[1000](82,18)(85,21)(88,24)
\qbezier[1000](88,24)(91,27)(94,30)
\qbezier[1000](94,30)(97,33)(100,36)
\qbezier[1000](100,36)(97,39)(94,42)
\qbezier[1000](94,42)(91,48)(88,54) \put(81,12){\makebox(0,0){\tiny
$Q_1$}} \put(93,24){\makebox(0,0){\tiny $Q_2$}}
\put(105,36){\makebox(0,0){\tiny $Q_3$}}
\put(3,-2){\dashbox{0.5}(18,23)} \put(4,4){\tiny right subtree of
$Q_2$}

\put(67,-2){\dashbox{0.5}(18,23)} \put(71,4){\tiny left subtree of
$Q_2$}
\end{picture}
\end{center}
\caption {The bijection $\rho$.} \label{fig1}
\end{figure}

Let $\mathcal{A}_n$ be the set of binary plane trees in
$\mathcal{T}_n$ where three vertices forming a $213$-pattern are
colored black. Let $\mathcal{B}_n$ be the set of all binary plane
trees in $\mathcal{T}_n$ where three vertices forming a
123-pattern are colored black. We will define a map
$\rho: \mathcal{A}_n \rightarrow \mathcal{B}_n$ as follows. Given a tree $T\in \mathcal{A}_n$ with $Q_2,Q_1,Q_3$ being the three
black vertices as a $213$-pattern, define $\rho(T)$ be the tree
obtained by changing the right subtree of $Q_2$ to be its left subtree.
See Figure~\ref{fig1} for an illustration.

In the tree $\rho(T)$, the relative positions of $Q_2$ and $Q_3$ keep the same, and
$Q_1$ is a left descendant of $Q_2$. Therefore, the three black points $Q_1Q_2Q_3$
form a $123$-pattern in $\rho(T)$, and $\rho(T)\in \mathcal{B}_n$. It is easy to
describe the converse and we omit here.

The first few values of $f_q(S_n(132,231))$ for $q$ of length $3$ are shown below.
\begin{center}
\begin{tabular}{|c|c|c|c|c|c|c||c|c|c|c|c|c|c|}
\hline $n$&$f_{123}$&$f_{132}$&$f_{213}$ &$f_{231}$&$f_{312}$  &$f_{321}$
      &$n$&$f_{123}$&$f_{132}$&$f_{213}$ &$f_{231}$&$f_{312}$  &$f_{321}$\\
\hline   $3$ &$1$      &$0$      &$1$       &$0$      &$1$&$1$  &$6$ &$160$      &$0$       &$160$     &$0$    &$160$     &$160$\\
\hline   $4$ &$8$      &$0$       &$8$       &$0$      &$8$       &$8$ &$7$ &$560$      &$0$       &$560$     &$0$    &$560$     &$560$\\
\hline   $5$ &$40$      &$0$       &$40$      &$0$     &$40$      &$40$  &$8$ &$1792$      &$0$       &$1792$     &$0$    &$1792$     &$1792$\\
\hline
\end{tabular}\label{2c}
\end{center}

\subsection{Pattern Count on $(132,312)$-Avoiding Permutations}
We begin with a correspondence between $(132,312)$-avoiding
permutations and compositions of $n$.
\begin{thm}
There is a bijection $\varphi_4$ between the sets $S_n(132,312)$ and
$\mathcal{C}_n$.
\end{thm}
\pf For $\sg\in S_n(132,312)$, let $\sg_{i_1},\sg_{i_2},\ldots, \sg_{i_k}$ be the
$k$ left to right maxima with $i_1<i_2<\cdots<i_k$, and thus,
$c=(i_2-i_1)+(i_3-i_2)+\cdots+(i_{k}-i_{k-1})+(n+1-i_{k})$ is a composition of $n$
from $i_1=1$. On the converse, let $n=c_k+c_{k-1}+\cdots+c_2+c_1\in \C_n$. For
$1\leq i \leq k$, if $c_i=1$ then set $\tau_i=n-i+1$; otherwise, set
$m_i=c_1+\cdots+c_{i-1}-i+2$ and $\tau_i=n-i+1,m_i+c_i-2,\ldots,m_i+1,m_i$. It is
easy to get $\sg=\tau_k,\tau_{k-1},\ldots,\tau_2,\tau_1\in S_n(132,312)$ as desired.
For example, if $9=3+1+2+3$, then $\sg=6\,5\,4\,7\,8\,3\,9\,2\,1$. \qed

\begin{prop}\label{pd}
For $n\geq3$, we have
\begin{align}\label{pd1}
f_{123}(n)=\sum_{k\geq 1 \atop c_1+c_2+\cdots+c_k=n} \sum_{i=1}^{k-2}c_i{k-i\choose
2}.
\end{align}
\end{prop}

\pf Given a permutation $\sg=\tau_k,\ldots,\tau_2,\tau_1$ in $S_n(132,312)$ whose
composition is given by $n=c_k+c_{k-1}+\cdots+c_2+c_1$, we see that the first
element in $\tau_i$ is larger than all the elements in $\tau_{j}$, whereas the other
elements in $\tau_i$ are smaller than the elements in $\tau_{j}$ for $i+1\leq j \leq
k$. The left to right maxima form an increasing subsequence and the other elements
form a decreasing subsequence. Thus we have $c_i$ choices to choose one element in
$\tau_{k-i+1}$ to play the role of $1$, and then ${k-i\choose 2}$ choices to choose
two left to right maxima in $\tau_j$ for $j<k-i+1$ to paly the role of $23$, summing
up all the number of $123$ patterns in subsequences $\tau_k,\ldots,\tau_2,\tau_1$
gives the formula \eqref{pd1}.\qed

\begin{thm}\label{td1}
For $n\geq 3$, in the set $S_n(132,312)$, we have
\begin{align}\label{fd1}
f_{123}(n)&=f_{321}(n)={n\choose 3}2^{n-3},\\\label{fd2}
f_{213}(n)&=f_{231}(n)={n\choose 3}2^{n-3}.
\end{align}
\end{thm}

\pf  From Lemma~\ref{oper}, we see that $\sg \in S_n(132,312)\Leftrightarrow \sg^{c}\in S_n(132,312)$, which follows $f_{123}(n)=f_{321}(n)$ and $f_{213}(n)=f_{231}(n)$ from $123^{c}=321$ and $213^{c}=231$, respectively.

To calculate $f_{123}(n)$, we have by Prop~\ref{pd}
\[f_{123}(n+1)=\sum_{k\geq 1,c_1=1 \atop c_1+c_2+\cdots+c_k=n+1}
\sum_{i=1}^{k-2}c_i{k-i\choose 2}+\sum_{k\geq 1,c_1\geq 2 \atop c_1+c_2+\cdots+c_k=n+1}
\sum_{i=1}^{k-2}c_i{k-i\choose 2}.
\]
If $c_1=1$, then $k\geq 2$, and
\begin{align*}
\sum_{k\geq 1,c_1=1 \atop c_1+c_2+\cdots+c_k=n+1}
\sum_{i=1}^{k-2}c_i{k-i\choose 2}&=\sum_{k-1\geq 1 \atop c_2+\cdots+c_k=n}
\sum_{i=2}^{k-2}c_i{k-i\choose 2}+\sum_{k-1\geq 1\atop c_2+\cdots+c_k=n}
{k-1\choose 2}\\
&=f_{123}(n)+d(n)/2.
\end{align*}
If $c_1\geq 2$, let $c'_1=c_1-1$, $c'_i=c_i$ for $2\leq i \leq k$, then $c'_1\geq 1$, and
\begin{align*}
\sum_{k\geq 1,c_1\geq 2 \atop c_1+c_2+\cdots+c_k=n+1}
\sum_{i=1}^{k-2}c_i{k-i\choose 2} &=\sum_{k\geq 1\atop c'_1+c'_2+\cdots+c'_k=n}
\sum_{i=1}^{k-2}c'_i{k-i\choose 2}
+\sum_{k\geq 1\atop c'_1+c'_2+\cdots+c'_k=n}
{k-1\choose 2}\\
&=f_{123}(n)+\sum_{k\geq 1\atop c'_1+c'_2+\cdots+c'_k=n}\left[{k\choose 2}+1-k\right]\\
&=f_{123}(n)+d(n)/2+2^{n-1}-c(n).
\end{align*}
Combining the above two cases, we get that
\[
f_{123}(n+1)=2f_{123}(n)+(n^2-n)2^{n-3},
\]
and the formula~\eqref{fd1} is derived by solving the recurrence with initial value $f_{123}(3)=1$. Formula~\eqref{fd2} is a direct computation of the equality $2f_{123}(n)+2f_{213}(n)={n\choose 3}2^{n-1}$.\qed

In the subsequent section, we could also give a combinatorial interpretation for $f_{231}(n)=f_{123}(n)$ by using binary plane trees.

For $\sg \in S_n(132,312)$, as in previous section, we could construct a binary
plane tree $T(\sg)$ on $n$ vertices such that each vertex that is a right descendant
of some vertex does not have a left descendant from the forbiddance of the pattern
$312$. Denote by $\mathscr{T}_n$ the set of such trees on $n$ vertices. Let $Q$ be
an occurrence of the pattern $231$ in $\sg \in S_n(132,312)$, and let $Q_2,Q_3,Q_1$
be the three vertices of $T(\sg)$ that correspond to $Q$, going left to right. Then,
$Q_2$ is a left descendant of $Q_3$, and there exists a lowest ascendant $x$ of
$Q_3$ or $x=Q_3$ so that $Q_1$ is a right descendant of $x$. Let $\mathscr{A}_n$ be
the set of binary plane trees in $\mathscr{T}_n$ in which three vertices forming a
$231$-pattern are colored black. Let $\mathscr{B}_n$ be the set of all binary plane
trees in $\mathscr{T}_n$ in which three vertices forming a 123-pattern are colored
black. It remains to construct a map $\varrho: \mathscr{A}_n \rightarrow
\mathscr{B}_n$.

\begin{figure}[h,t]
\setlength{\unitlength}{0.7mm}%
\begin{center}
\begin{picture}(50,55)
\put(0,0){\circle*{2}}\put(12,0){\circle*{2}}
\put(6,12){\circle*{3.5}}\put(9,18){\circle*{2}}
\put(12,24){\circle*{3.5}}
\put(15,30){\circle*{2}}
\put(18,36){\circle*{2}}
\put(21,30){\circle*{2}}
\put(24,24){\circle*{3.5}}
\put(30,12){\circle*{2}}
\put(21,42){\circle*{2}}
\put(24,48){\circle*{2}}\put(27,42){\circle*{2}}
\qbezier[1000](0,0)(3,6)(6,12)\qbezier[1000](12,0)(9,6)(6,12)
\qbezier[1000](6,12)(7.5,15)(9,18)
\qbezier[1000](9,18)(10.5,21)(12,24)
\qbezier[1000](12,24)(13.5,27)(15,30)
\qbezier[1000](15,30)(16.5,33)(18,36)
\qbezier[1000](18,36)(19.5,33)(21,30)
\qbezier[1000](21,30)(22.5,27)(24,24)
\qbezier[1000](24,24)(27,18)(30,12)
\qbezier[1000](18,36)(19.5,39)(21,42)
\qbezier[1000](21,42)(22.5,45)(24,48)
\qbezier[1000](24,48)(25.5,45)(27,42)
\put(2,12){\makebox(0,0){\tiny $Q_2$}}
\put(8,24){\makebox(0,0){\tiny $Q_3$}}
\put(29,24){\makebox(0,0){\tiny $Q_1$}}
\put(16,35){\makebox(0,0){\tiny $x$}}
\put(23,42){\makebox(0,0){\tiny $y$}}
\put(45,18){\makebox(0,0){$\Rightarrow$}}

\put(60,0){\circle*{2}}\put(72,0){\circle*{2}}
\put(66,12){\circle*{3.5}}\put(69,18){\circle*{2}}
\put(72,24){\circle*{3.5}}
\put(75,30){\circle*{2}}
\put(78,36){\circle*{2}}\put(81,42){\circle*{3.5}}
\put(81,30){\circle*{2}}
\put(87,30){\circle*{2}}
\put(84,48){\circle*{2}}\put(87,54){\circle*{2}}
\put(90,48){\circle*{2}}

\qbezier[1000](60,0)(63,6)(66,12)\qbezier[1000](72,0)(69,6)(66,12)
\qbezier[1000](66,12)(67.5,15)(69,18)
\qbezier[1000](69,18)(70.5,21)(72,24)
\qbezier[1000](72,24)(73.5,27)(75,30)
\qbezier[1000](75,30)(76.5,33)(78,36)
\qbezier[1000](78,36)(79.5,33)(81,30)

\qbezier[1000](78,36)(79.5,39)(81,42)
\qbezier[1000](81,42)(84,36)(87,30)
\qbezier[1000](81,42)(82.5,45)(84,48)
\qbezier[1000](84,48)(85.5,51)(87,54)
\qbezier[1000](87,54)(88.5,51)(90,48)

\put(62,12){\makebox(0,0){\tiny $Q_2$}}
\put(68,24){\makebox(0,0){\tiny $Q_3$}}
\put(77,42){\makebox(0,0){\tiny $Q_1$}}
\put(76,35){\makebox(0,0){\tiny $x$}}
\put(86.5,48){\makebox(0,0){\tiny $y$}}

\put(19,39){\dashbox{0.5}(10,13)}
\put(35,50){\footnotesize $T_u$}\put(26,48){$\nearrow$}

\put(82,45){\dashbox{0.5}(10,13)}
\put(97,50){\footnotesize $T_u$}
\put(90,50){$\rightarrow$}

\put(-2,-2){\dashbox{0.5}(24,40)}
\put(27,4){\footnotesize $T_d$}
\put(20,10){$\searrow$}

\put(58,-2){\dashbox{0.5}(24,40)}
\put(90,4){\footnotesize $T_d$}
\put(80,4){$\rightarrow$}

\end{picture}
\end{center}
\caption {The bijection $\varrho$.} \label{fig2}
\end{figure}
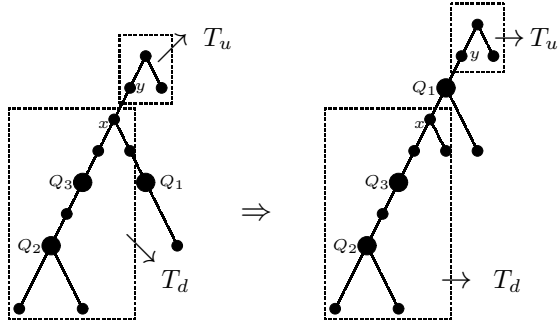

Given a tree $T\in \mathscr{A}_n$ with $Q_2,Q_3,Q_1$ being the three black vertices as a $231$-pattern, denoted by $y$ the vertex that is the parent of $x$ if it exists. We can see that $x$ is the left child of $y$ from $T\in \mathscr{A}_n$. Let $T^u:=T-T_x$ be the tree obtained from $T$ by deleting the subtree $T_x$, and
$T^d:=T_x-T_{Q_1}$ be the tree obtained from $T_x$ by deleting the
subtree $T_{Q_1}$. Now we can define $\varrho(T)$ be the tree obtained from $T$ by first interchanging $Q_1$ as the left child of $y$, then adjoining the subtree $T^d$ as the left subtree of the vertex $Q_1$ and keeping all three black vertices the same if $y$ exits. Otherwise, we adjoin the subtree $T^d$ as the left subtree of the vertex $Q_1$ directly. An illustration is given in Figure~\ref{fig2}.

In the tree $\varrho(T)$, the relative positions of $Q_2$ and $Q_3$
are unchanged, and $Q_3$ is a left descendant of $Q_1$, thus the
three black points $Q_2Q_3Q_1$ form a $123$-pattern in $\varrho(T)$, and $\varrho(T)\in \mathscr{B}_n$. It is easy to describe the converse map and we omit here.

The first few values of $f_q(S_n(132,312))$ for $q$ of length $3$ are shown below.
\begin{center}
\begin{tabular}{|c|c|c|c|c|c|c||c|c|c|c|c|c|c|}
\hline $n$ &$f_{123}$&$f_{132}$&$f_{213}$ &$f_{231}$&$f_{312}$  &$f_{321}$ &$n$ &$f_{123}$&$f_{132}$&$f_{213}$ &$f_{231}$&$f_{312}$  &$f_{321}$\\
\hline   $3$ &$1$      &$0$      &$1$       &$1$        &$0$      &$1$& $6$&$160$      &$0$       &$160$     &$160$   &$0$         &$160$\\
\hline   $4$ &$8$      &$0$       &$8$       &$8$     &$0$             &$8$  &$7$ &$560$      &$0$       &$560$     &$560$   &$0$         &$560$\\
\hline   $5$ &$40$      &$0$       &$40$      &$40$    &$0$           &$40$  &$8$ &$1792$      &$0$       &$1792$     &$1792$   &$0$         &$1792$\\
\hline
\end{tabular}\label{2d}
\end{center}

\subsection{Pattern Count on $(132,321)$-Avoiding Permutations}
We begin with a correspondence from Simion and Schmidt \cite{Sim} as follows:
\begin{lem}\label{2le}
There is a bijection $\varphi_5$ between the set
$S_n(132,321)\backslash \{\text{identity}\}$ and the set of
$2$-element subsets of $[n]$.
\end{lem}
\pf For a permutation $\sg \in S_n(132,321)\backslash
\{\text{identity}\}$, suppose that $\sg_k=m$ ($k<m$), and then we define
$\varphi_5(\sg)=\{k,m\}$. On the converse, given two elements $1\leq
k<m \leq n$, set $\tau_1=m-k+1,m-k+2,\ldots,m-1,m$,
$\tau_2=1,2,\ldots,m-k$ and $\tau_3=m+1,m+2,\ldots,n-1,n$. Then
define $\sg=\varphi_5^{-1}(k,m)=\tau_1,\tau_2,\tau_3$. For example, if $k=4,m=6$, then $\sg=3\,4\,5\,6\,1\,2\,7\,8$. \qed

From the above lemma, we have

\begin{prop}\label{pe}
For $n\geq3$,
\begin{align}
f_{213}(n)&=\sum_{1\leq k<m \leq n}k(m-k)(n-m)\label{pe1},\\
f_{312}(n)&=\sum_{1\leq k<m \leq n}k{m-k\choose 2}.~~~~~~\label{pe2}
\end{align}
\end{prop}
\pf Given a permutation $\sg=\tau_1,\tau_2,\tau_3$ with
$\varphi_5(\sg)=\{k,m\}$ in the set $S_n(132,321)$, we see that the
elements in $\tau_1$, $\tau_2$ and $\tau_3$ are increasing, and
$\tau_2<\tau_1<\tau_3$. We have $k$ choices to select one
element in $\tau_{1}$ to play the role of $2$, and then have $m-k$ choices to
choose one element in $\tau_{2}$ to play the role of $1$, and $n-m$
choices to choose one element in $\tau_{3}$ to play the role of $3$.
Summing up all possible $k$ and $m$ gives the formula \eqref{pe1}.

For the pattern $312$, we have $k$ choices in the subsequence
$\tau_1$ to choose one element to play the role of $3$, after this
we have ${m-k\choose 2}$ choices in the subsequence $\tau_2$ to
choose two elements to play the role of $23$. Summing up all
$k$ and $m$ gives the formula \eqref{pe2}. \qed

Next, we exhibit some useful formulae for calculating $f_{213}(n)$ and $f_{312}(n)$ as follows:
\begin{lem}\label{le}
For $n\geq 2$,
\begin{align*}
\sum_{k=1}^{n-1}k&=\frac{n(n-1)}{2};~~~~~
\sum_{k=1}^{n-1}k^2=\frac{n(n-1)(2n-1)}{6};\\
\sum_{k=1}^{n-1}k^3&=\frac{n^2(n-1)^2}{4};~~~
\sum_{k=1}^{n-1}k^4=\frac{n(n-1)(2n-1)(3n^2-3n+1)}{30}.
\end{align*}
\end{lem}

Based on previous analysis, we are ready to give the exact formulae for the pattern count in $S_n(132,321)$.
\begin{thm}\label{te1}
For $n\geq 3$, in the set $S_n(132,321)$, we have
\begin{align}
f_{213}(n)&=f_{231}(n)=f_{312}(n)={n+2\choose 5},~~~~\\
f_{123}(n)&=n(7n^4-40n^3+85n^2-80n+28)/120.
\end{align}
\end{thm}
\pf From Lemma~\ref{oper}, we see that $\sg \in
S_n(132,321)\Leftrightarrow \sg^{-1}\in S_n(132,321)$, it implies that $f_{312}(n)=f_{231}(n)$ since $312^{-1}=231$. By using Prop~\ref{pe}, we have
\begin{align*}
f_{312}(n)&=\sum_{1\leq k<m \leq n}k{m-k\choose 2}\\
&=\sum_{k=1}^{n-1}k\sum_{m=k+1}^n{m-k\choose 2}=\sum_{k=1}^{n-1}k{n-k+1\choose 3}\\
&=\sum_{k=1}^{n-1}\left[(n^3-n)k+(1-3n^2)k^2+2nk^3-k^4\right]={n+2\choose 5},
\end{align*}
where the last equality holds from Lemma~\ref{le}. By using Prop~\ref{pe} again, we have
\begin{align*}
f_{213}(n)&=\sum_{1\leq k<m \leq n}k(m-k)(n-m)=\sum_{k=1}^{n-1}\sum_{m=k+1}^nk(m-k)(n-m)\\
&=\sum_{k=1}^{n-1}\sum_{m'=1}^{n-k}km'(n-m'-k)=\sum_{k=1}^{n-1}k(n-k)\sum_{m'=1}^{n-k}m'-\sum_{k=1}^{n-1}k\sum_{m'=1}^{n-k}m'^2\\
&=\sum_{k=1}^{n-1}\left[ \left(\frac{n^3}{6}-\frac{n}{6}\right)k+\left(\frac{1}{6}-\frac{n^2}{2}\right)k^2+ \frac{n}{2}k^3-\frac{1}{6}k^4\right]={n+2\choose 5},
\end{align*}
where the last equality holds by simple calculation from Lemma~\ref{le}. We complete the proof by employing the relation $2f_{231}(n)+f_{213}(n)+f_{123}(n)={n\choose 3}\left[{n\choose 2}+1\right]$.\qed

We also notice that the equality $f_{213}(n)=f_{231}(n)$ can be proved by using B\'{o}na's bijection \cite{Bona} directly on the set of binary plane trees on $n$ vertices such that the vertex which is a right descendant of some node has no right descendants.

The first few values of $f_q(S_n(132,321))$ for $q$ of length $3$ are shown below.
\begin{center}
\begin{tabular}{|c|c|c|c|c|c|c||c|c|c|c|c|c|c|}
\hline  $n$&$f_{123}$&$f_{132}$&$f_{213}$ &$f_{231}$&$f_{312}$  &$f_{321}$
&$n$&$f_{123}$&$f_{132}$&$f_{213}$ &$f_{231}$&$f_{312}$  &$f_{321}$\\
\hline   $3$ &$1$      &$0$      &$1$     &$1$&$1$  &$0$      &$6$ &$152$      &$0$       &$56$   &$56$ &$56$   &$0$ \\
\hline   $4$ &$10$      &$0$       &$6$   &$6$ &$6$     &$0$   &$7$ &$392$      &$0$       &$126$  &$126$ &$126$   &$0$   \\
\hline   $5$ &$47$      &$0$       &$21$   &$21$  &$21$     &$0$    &$8$ &$868$      &$0$       &$252$  &$252$ &$252$   &$0$ \\
\hline
\end{tabular}\label{2e}
\end{center}

\section{Triply Restricted Permutations}
In this section, we study the pattern count in the simultaneous avoidance of any three patterns of length $3$. Based on Lemma~\ref{oper}, Simion and Schmidt \cite{Sim} showed that the pairs of patterns among the total ${6\choose 3}=20$ cases fall into the following $6$ classes:
\begin{prop}
For every symmetric group $S_n$,\\{\small
(1) $|S_n(123,132,213)|=|S_n(231,312,321)|=F_{n+1}$;\\
(2) $|S_n(123,132,231)|=|S_n(123,213,312)|=|S_n(132,231,321)|=|S_n(213,312,321)|=n$; \\
(3) $|S_n(132,213,231)|=|S_n(132,213,312)|=|S_n(132,231,312)|=|S_n(213,231,312)|=n$; \\
(4) $|S_n(123,132,312)|=|S_n(123,213,231)|=|S_n(132,312,321)|=|S_n(213,231,321)|=n$; \\
(5) $|S_n(123,231,312)|=|S_n(132,213,321)|=n$;\\
(6) $|S_n(R)|=0$ for all $R\supset \{123,321\}$ if $n\geq 5$,\\
where $F_{n}$ is the Fibonacci number given by $F_0=0,F_1=1$ and $F_{n}=F_{n-1}+F_{n-2}$ for $n\geq 2$.}
\end{prop}

\subsection{Pattern Count on $(123,132,213)$-Avoiding Permutations}

It is known that $F_{n+1}$ counts the number of $0$-$1$ sequences of length $n-1$ in which there are no consecutive ones, see \cite{Com}, and we call such a sequence a Fibonacci binary word for convenience. Let $B_n$ denote the set of all Fibonacci binary words of length $n$. Simion and Schmidt \cite{Sim} showed that
\begin{lem}[\cite{Sim}]
There is a bijection $\psi_1$ between $S_n(123,132,213)$ and $B_{n-1}$.
\end{lem}
\pf Let $w=w_1w_2\cdots w_{n-1} \in B_{n-1}$, and the corresponding permutation $\sg$ is determined as follows: For $1\leq i \leq n-1$, let $X_i=[n]-\{\sg_1,\ldots,\sg_{i-1}\}$, and then set
\begin{subnumcases}
{\sg_i=}
\text{largest element in}~ X_i, &if $s_i=0$, \nonumber\\
\text{second largest element in}~ X_i, &if $s_i=1$.\nonumber
\end{subnumcases}
Finally, $\sg_n$ is the unique element in $X_n$. For example, if $w=01001010$, then $\psi_1(w)=9\,7\,8\,6\,4\,5\,2\,3\,1$.\qed

Given a word $w=w_1w_2\cdots w_{n} \in B_{n}$, we call $i$ ($1\leq
i<n$) an ascent of $w$ if $w_i<w_{i+1}$, and denote by $\asc(w)=\{ i|
w_i<w_{i+1} \}$ and $\maj(w)=\sum\limits_{i \in \asc(w)}i$.
\begin{prop}\label{paa}
The total number of occurrences of the pattern $312$ in $S_n(123,132,213)$ is given by
\begin{align}\label{paa1}
f_{312}(n)=\sum_{w \in B_{n-1}}\maj(w).
\end{align}
\end{prop}

\pf Suppose $\sg \in S_n(123,132,213)$ and $\psi_1(\sg)=w_1w_2\cdots
w_{n-1}$. If $k$ is an ascent of $w$, then $w_kw_{k+1}=01$. We have $\sg_k>\sg_{j}$ for all $j>k$ since $\sg_{k}$ is the largest element in $X_{k}$. On the other hand, there exists a unique $l>k+1$ such that $\sg_{l}>\sg_{k+1}$ since $\sg_{k+1}$ is the second largest element in $X_{k+1}$. From the bijection $\psi_1$, we see that
for all $i\in [n-1]$, there is at most one $j>i$ such that
$\sg_j>\sg_i$; this implies that $\sg_i>\sg_{k+1}$ for all $i<k$.
Thus we find that $\sg_i \sg_{k+1} \sg_{l}$ forms a $312$-pattern
for all $i\leq k$, that is the ascent $k$ will produce $k$'s
$312$-pattern in which $\sg_{k+1}$ plays the role of $1$. Summing up
all the ascents, we derive that there are total $\maj(w)$ such
patterns in $\sg$.  \qed

\begin{thm}\label{Tta1}
For $n\geq 3$, in the set $S_n(123,132,213)$, we have
\begin{align}
\sum_{n\geq 3}f_{231}(n)x^n&=\sum_{n\geq 3}f_{312}(n)x^n=\frac{x^3(1+2x)}{(1-x-x^2)^3},\label{fTta1}\\
\sum_{n\geq 3}f_{321}(n)x^n&= \frac{x^3(1+6x+12x^2+8x^3)} {(1-x-x^2)^4}.\label{Tca2}
\end{align}
\end{thm}

\pf From Lemma~\ref{oper}, we have $f_{231}(n)=f_{312}(n)$ since $\sg \in
S_n(123,132,213)\Leftrightarrow \sg^{-1}\in S_n(123,132,213)$ and $231^{-1}=312$. By Prop~\ref{paa}, we write
\[
\sum_{n\geq 3}f_{312}(n)x^n=\sum_{n\geq 3}x^n\sum_{\sg \in B_{n-1}}\maj(w) =x\sum_{n\geq 3}\sum_{w \in B_{n-1}}\maj(w)x^{n-1}=xu(x),
\]
where $u(x)=\sum\limits_{n\geq 2}\sum\limits_{w \in B_n}\maj(w)x^n$.

Let $M_n(q)=\sum\limits_{w \in B_n}  q^{\maj(w)}$ and $M(x,q)=\sum\limits_{n\geq
2}M_n(q)x^n$. Then $u(x)=\frac{\partial M(x,q)}{\partial
q}\mid_{q=1}$. Given a word $w=w_1w_2\cdots w_n \in B_n$, if
$w_n=0$, then $\maj(w)=\maj(w_1w_2\cdots w_{n-1})$; otherwise, $w_{n-1}w_n=01$ and $\maj(w)=\maj(w_1w_2\cdots w_{n-2})+(n-1)$. Hence, we have
\begin{align}
M_n(q)=M_{n-1}(q)+q^{n-1}M_{n-2}(q) \text{ for }n\geq 4,
\end{align}
with $M_2(q)=2+q$ and $M_3(q)=2+q+2q^2$.
Multiplying the recursion by $x^n$ and summing over $n\geq 4$ yields that
\[ M(x,q)-(2+q)x^2-(2+q+2q^2)x^3=x\left[M(x,q)-(2+q)x^2\right]+qx^2M(xq,q).\]
Therefore \[(1-x)M(x,q)=qx^2M(xq,q)+(2+q)x^2+2q^2x^3.\]
Differentiate both sides with respect to $q$, we get
\begin{align}\label{ftaf}
(1-x)\frac{\partial M(x,q)}{\partial q}=x^2\left[M(xq,q)+q\frac{\partial M(xq,q)}{\partial q} \right]+x^2+4qx^3.
\end{align}
Setting $q=1$ in equation \eqref{ftaf} reads that
\[
(1-x)u(x)=x^2\left[M(x,1)+\frac{\partial M(xq,q)}{\partial q}\mid_{q=1}\right]
+x^2+4x^3,
\]
Employing the well-known generating  $\sum_{n\geq 0}F_nx^n=\frac{x}{1-x-x^2}$, we have
\[
M(x,1)=\sum_{n\geq 2}|B_{n}|x^n=\sum_{n\geq 2}F_{n+2}x^n= \frac{x^2(3+2x)}{1-x-x^2}.
\]
Further,
\begin{align*} \frac{\partial M(xq,q)}{\partial
q}\mid_{q=1}&=\left(\sum_{n\geq 2}\sum_{w \in B_n}  q^{n+\maj(w)}x^n
\right)|_{q=1}\\
&=\sum_{n\geq 2}x^n\sum_{w \in B_n}(n+\maj(w))=\sum_{n\geq
2}nF_{n+2}x^n+u(x).
\end{align*}
From the generating function of $F_{n+2}$, we obtain
\[
\sum_{n\geq 2}nF_{n+2}x^n=x\left(\frac{x^2(3+2x)}{1-x-x^2}\right)' =\frac{x^2(6+3x-4x^2-2x^3)}{(1-x-x^2)^2},
\]
which yields that
\[(1-x)u(x)=x^2\left[\frac{x^2(3+2x)}{1-x-x^2}+\frac{x^2(6+3x-4x^2-2x^3)}{(1-x-x^2)^2}+u(x)\right]
+x^2+4x^3.\]
Therefore,
\[
u(x)={x^2(1+2x)}/{(1-x-x^2)^3},
\]
which gives the generating function for $f_{312}(n)$ as shown in formula~\eqref{fTta1}.

For formula~\eqref{Tca2}, we first have
\begin{align}\label{Tta2f}
\sum_{n\geq 3}f_{321}(n)x^n=\sum_{n\geq 3}{n\choose 3}F_{n+1}x^n-
2\sum_{n\geq 3}f_{312}(n)x^n
\end{align}
from $2f_{312}(n)+f_{321}(n)={n\choose 3}F_{n+1}$. Using the fact that $\sum_{n\geq 0}F_nx^n=\frac{x}{1-x-x^2}$, we get
\begin{align*}
&\sum_{n\geq 3}F_{n+1}x^n=\frac{1}{x}\left(\frac{x}{1-x-x^2}-x-x^2-2x^3\right)=\frac{x^3(3+2x)}{1-x-x^2},\\
&\sum_{n\geq 3}{n\choose 3}F_{n+1}x^n=\frac{x^3}{6}\left(\sum_{n\geq 3}F_{n+1}x^n\right)^{'''}=\frac{x^3(3+8x+6x^2+4x^3)}{(1-x-x^2)^4}.
\end{align*}
Formula~\eqref{Tca2} follows by substituting the generating functions of ${n\choose 3}F_{n+1}$ and $f_{312}(n)$ into \eqref{Tta2f}, and we complete the proof.\qed

The first few values of $f_q(S_n(123,132,213))$ for $q$ of length $3$ are shown below.
\begin{center}
\begin{tabular}{|c|c|c|c|c|c|c||c|c|c|c|c|c|c|}
\hline  $n$&$f_{123}$&$f_{132}$&$f_{213}$ &$f_{231}$&$f_{312}$  &$f_{321}$
& $n$&$f_{123}$&$f_{132}$&$f_{213}$ &$f_{231}$&$f_{312}$  &$f_{321}$\\
\hline   $3$ &$0$      &$0$      &$0$       &$1$      &$1$      &$1$ & $6$ &$0$      &$0$       &$0$     &$40$    &$40$     &$180$\\
\hline   $4$ &$0$      &$0$       &$0$       &$5$      &$5$       &$10$ &$7$ &$0$      &$0$       &$0$     &$95$    &$95$     &$545$\\
\hline   $5$ &$0$      &$0$       &$0$      &$15$     &$15$      &$50$ &$8$ &$0$      &$0$       &$0$     &$213$    &$213$     &$1478$\\
\hline
\end{tabular}\label{3a}
\end{center}

\subsection{Pattern Count on Other Triple Avoiding Permutations}
We begin with the pattern count on $(123,132,231)$-avoiding permutations.
\begin{thm}\label{Ttb1}
For $n\geq 3$, in the set $S_n(123,132,231)$, we have
\begin{align}
f_{213}(n)&=f_{312}(n)={n\choose 3},\\
f_{321}(n)&=(n-2){n\choose 3}.
\end{align}
\end{thm}

\pf We first give the following structure from Simion and Schmidt \cite{Sim}
\begin{equation}\label{eq:temp3}
\sg\in S_n(123,132,231)\Leftrightarrow \sg=n,n-1,\ldots,k+1,k-1,k-2,\ldots,1,k \text{ for some } k.
\end{equation}
Based on such structure, we can show $f_{213}(n)=f_{312}(n)$ by a direct bijection.  Let $q=abc$ be a $213$-pattern of a permutation $\sg \in S_n(123,132,231)$. From $b<c$ and the fact that $\sg\in S_n(123,132,231)$ has only one ascent at position $n-1$, we have $\sg(n)=c$, and thus $\sg=n,n-1,\ldots,c+1,c-1,c-2,\ldots,2,1,c$ from the structure~\eqref{eq:temp3}. Let $q'=cba$ (a $312$-pattern) and we set $\sg'=n,n-1,\ldots,a+1,a-1,a-2,\ldots,2,1,a$ as the desired permutation. For example, if $n=7$ and $q=326$, then $\sg=7\,5\,4\,3\,2\,1\,6$. Moreover, $q'=623$ and $\sg'=7\,6\,5\,4\,2\,1\,3$.

To calculate $f_{312}(n)$, we suppose $\sg=n,n-1,\ldots,k+1,k-1,k-2, \ldots,2,1,k$ for some $k$ from structure~\eqref{eq:temp3}. We can construct a $312$-pattern as follows: Choose
one element from the first $n-k$ elements to play the role of $3$,
then choose one element from the next $k-1$ elements to play the
role of $1$, and the last element plays the role of $2$. Thus,
summing up all possible $k$, we have
\[
f_{312}(n)=\sum_{k=1}^n
(n-k)(k-1)=-n^2+(n+1)\sum_{k=1}^nk-\sum_{k=1}^nk^2=\frac{n(n-1)(n-2)}{6}={n\choose
3}.
\]
The formula for $f_{321}(n)$ follows by using $f_{213}(n)+f_{312}(n)+f_{321}(n)=n{n\choose 3}$.\qed

The first few values of $f_q(S_n(123,132,231))$ for $q$ of length $3$ are shown below.
\begin{center}
\begin{tabular}{|c|c|c|c|c|c|c||c|c|c|c|c|c|c|}
\hline  $n$&$f_{123}$&$f_{132}$&$f_{213}$ &$f_{231}$&$f_{312}$  &$f_{321}$ &$n$&$f_{123}$&$f_{132}$&$f_{213}$ &$f_{231}$&$f_{312}$  &$f_{321}$\\
\hline   $3$ &$0$      &$0$      &$1$&$0$             &$1$      &$1$ &$6$ &$0$      &$0$       &$20$&$0$         &$20$     &$80$\\
\hline   $4$ &$0$      &$0$      &$4$ &$0$             &$4$       &$8$&$7$ &$0$      &$0$       &$35$   &$0$         &$35$     &$175$\\
\hline   $5$ &$0$      &$0$      &$10$ &$0$           &$10$      &$30$ &$8$ &$0$      &$0$       &$56$   &$0$         &$56$     &$336$\\
\hline
\end{tabular}\label{3b}
\end{center}

For $(132,213,231)$-avoiding permutations, we have
\begin{thm}\label{Ttb1}
For $n\geq 3$, in the set $S_n(132,213,231)$,
\begin{align}
f_{123}(n)&=f_{312}(n)={n+1\choose 4},\\
f_{321}(n)&=\frac{n(n-2)(n-1)^2}{12}.
\end{align}
\end{thm}
\pf We begin with the following structure by Simion and Schmidt \cite{Sim}
\begin{equation}\label{eq:tmp5}
\sg\in S_n(132,213,231)\Leftrightarrow
 \sg=n,n-1,\ldots,k+1,1,2,\ldots,k-1,k \text{~for some~} k.
\end{equation}
Based on such structure, we first prove $f_{123}(n)=f_{312}(n)$. For each
\[
\sg=n,n-1,\ldots,k+1,1,\ldots, \underline{a},a+1,\ldots, \underline{b},b+1,\ldots,c-1,\underline{c},c+1,\ldots,k-1,k
\]
with $abc$ as a $123$-pattern, we set
\[
\sg'=n,n-1,\ldots,\underline{n-k+c},\ldots,c,1,2,\ldots, \underline{a},a+1,\ldots,\underline{b},b+1,\ldots,c-1,
\]
and it is easy to check that $n-k+c,a,b$ is a $312$-pattern of $\sg'$. For example, if $\sg=9\,8\,7\,1\,\underline{2}\, \underline{3}\,4\, \underline{5}\,6$ then $\sg'=9\,\underline{8}\,7\,6\,5\,1\, \underline{2}\, \underline{3}\,4$.

To calculate $f_{123}(n)$, we suppose $\sg=n,n-1,\ldots,k+1,1,2,\ldots,k-1,k$ for some $k$ from structure~\eqref{eq:tmp5}. Then,  a $123$-pattern can be obtained by choosing three elements from the last $k$ elements to play the role of $123$. Thus, summing up all possible $k$ gives
\[
f_{123}(n)=\sum_{k=1}^n {k \choose 3}={n+1\choose 4}.
\]
We complete the proof by using $f_{123}(n)+f_{312}(n)+f_{321}(n)=n{n\choose 3}$.\qed

The first few values of $f_q(S_n(132,213,231))$ for $q$ of length $3$ are shown below.
\begin{center}
\begin{tabular}{|c|c|c|c|c|c|c||c|c|c|c|c|c|c|}
\hline $n$ &$f_{123}$&$f_{132}$&$f_{213}$ &$f_{231}$&$f_{312}$  &$f_{321}$ &$n$ &$f_{123}$&$f_{132}$&$f_{213}$ &$f_{231}$&$f_{312}$  &$f_{321}$\\
\hline   $3$ &$1$      &$0$      &$0$       &$0$      &$1$      &$1$ &$6$ &$35$      &$0$       &$0$     &$0$    &$35$     &$50$\\
\hline   $4$ &$5$      &$0$       &$0$       &$0$      &$5$       &$6$  & $7$ &$70$      &$0$       &$0$     &$0$    &$70$     &$105$\\
\hline   $5$ &$15$      &$0$       &$0$      &$0$     &$15$      &$20$ & $8$ &$126$      &$0$       &$0$     &$0$    &$126$     &$196$\\
\hline
\end{tabular}\label{3c}
\end{center}

\

For $(123,132,312)$-avoiding permutations, we have
\begin{thm}\label{Ttb1}
For $n\geq 3$, in the set $S_n(123,132,312)$,
\begin{align}
f_{213}(n)&=f_{231}(n)={n\choose 3},\\
f_{321}(n)&=(n-2){n\choose 3}.
\end{align}
\end{thm}
\proof We begin with the following structure from Simion and Schmidt
\cite{Sim}
\begin{equation}\label{eq:tmp4}
\sg\in S_n(132,213,231)\Leftrightarrow \sg=n-1,n-2,\ldots,k+1,n,k,k-1,\ldots,1 \text{ for some } k.
\end{equation}
Based on such structure, we prove $f_{213}(n)=f_{231}(n)$  by a direct correspondence. Let $\sg=n-1,\ldots,\underline{a},a+1,\ldots,\underline{b},b+1,\ldots,k+1,\underline{n}, k,k-1, \ldots,2,1$ with $abn$ as a $213$-pattern. Then, we set
\[\sg'=n-1,\ldots,\underline{n-a+b},\ldots,n-a+k+1,\underline{n},n-a+k,n-a+k-1,\ldots, \underline{n-a},\ldots,2,1,\]
where $n-a+b,n,n-a$ is a $231$-pattern of $\sg'$. For example, if $\sg=8\,\underline{7}\,6\,\underline{5}\,4\,\underline{9}\,3\,2\,1$, then $\sg'=\sg=8\,\underline{7}\,6\,\underline{9}\,5\,4\,3\,\underline{2}\,1$.

To calculate $f_{213}(n)$, we suppose that $\sg=n-1,n-2,\ldots,k+1,n,k,k-1,\ldots,2,1$ for some $k$. Then, a $213$-pattern can be obtained by choosing two elements from the first $n-k-1$ elements to play the role of $21$,
and let $n$ play the role of $3$. Thus, summing up all possible
$k$, we have
\[
f_{213}(n)=\sum_{k=0}^{n-1} {n-k-1 \choose 2}={n\choose 3}.
\]
We complete the proof by using $f_{213}(n)+f_{231}(n)+f_{321}(n)=n{n\choose 3}$.\qed

The first few values of $f_q(S_n(123,132,312))$ for $q$ of length
$3$ are shown below.
\begin{center}
\begin{tabular}{|c|c|c|c|c|c|c||c|c|c|c|c|c|c|}
\hline   $n$&$f_{123}$&$f_{132}$&$f_{213}$ &$f_{231}$&$f_{312}$  &$f_{321}$ &$n$&$f_{123}$&$f_{132}$&$f_{213}$ &$f_{231}$&$f_{312}$  &$f_{321}$\\
\hline   $3$ &$0$      &$0$      &$1$       &$1$      &$0$      &$1$    &$6$ &$0$      &$0$       &$20$     &$20$    &$0$     &$80$\\
\hline   $4$ &$0$      &$0$       &$4$       &$4$      &$0$       &$8$  &$7$ &$0$      &$0$       &$35$     &$35$    &$0$     &$175$\\
\hline   $5$ &$0$      &$0$       &$10$      &$10$     &$0$      &$30$  &$8$ &$0$      &$0$       &$56$     &$56$    &$0$     &$336$\\
\hline
\end{tabular}\label{3d}
\end{center}

Finally, we study the pattern count on $(123,231,312)$-avoiding permutations.
\begin{thm}
For $n\geq 3$, in the set $S_n(123,231,312)$, we have
\begin{align}
f_{132}(n)&=f_{213}(n)={n+1\choose 4},\label{Tte1}\\
f_{321}(n)&=\frac{n(n-2)(n-1)^2}{12}.
\end{align}
\end{thm}
\pf From Lemma~\ref{oper}, we see that
\[
\sg \in S_n(123,231,312)\Leftrightarrow \sg^{r}\in S_n(321,132,213)\Leftrightarrow (\sg^{r})^c\in S_n(123,231,312).
\]
We have $f_{213}(n)=f_{132}(n)$ from $(213^{r})^c=312^c=132$.
The following structure of the set $S_n(123,231,312)$ is given by Simion and Schmidt \cite{Sim}
\begin{equation}\label{s3e}
\sg\in S_n(132,213,231)\Leftrightarrow
 \sg=k-1,k-2,\ldots,3,2,1,n,n-1\ldots,k \text{~for some~} k.
\end{equation}
Suppose that $\sg=k-1,k-2,\ldots,3,2,1,n,n-1\ldots,k$ for some $k$,
then a $213$-pattern can be obtained as follows: Choose two elements
from the first $k-1$ elements to play the role of $21$, and choose
one element from the last $n-k+1$ elements to play the role of $3$.
Thus, summing up all possible $k$, we have
\begin{align*}
f_{213}(n)&=\sum_{k=1}^n {k-1 \choose 2}(n-k+1)=\sum_{k=0}^{n-1}{k\choose 2}(n-k)={n+1\choose 4}.
\end{align*}
The formula for $f_{321}(n)$ is obtained by the relation $2f_{213}(n)+f_{321}(n)=n{n\choose 3}$.  \qed

The first few values of $f_q(S_n(123,231,312))$ for $q$ of length
$3$ are shown below.
\begin{center}
\begin{tabular}{|c|c|c|c|c|c|c||c|c|c|c|c|c|c|}
\hline  $n$&$f_{123}$&$f_{132}$&$f_{213}$ &$f_{231}$&$f_{312}$  &$f_{321}$ &$n$&$f_{123}$&$f_{132}$&$f_{213}$ &$f_{231}$&$f_{312}$  &$f_{321}$\\
\hline   $3$ &$0$       &$1$      &$1$  &$0$      &$0$          &$1$  &$6$ &$0$     &$35$    &$35$   &$0$       &$0$        &$50$\\
\hline   $4$ &$0$     &$5$      &$5$ &$0$       &$0$              &$6$  &$7$ &$0$     &$70$    &$70$    &$0$       &$0$       &$105$\\
\hline   $5$ &$0$     &$15$     &$15$    &$0$       &$0$         &$20$  &$8$ &$0$     &$126$    &$126$    &$0$       &$0$       &$196$\\
\hline
\end{tabular}\label{3e}
\end{center}


\end{document}